# On a SIS model with impulsive effects

M. De la Sen*, A. Ibeas ** and S. Alonso-Quesada*

*Department of Electricity and Electronics. Faculty of Science and Technology

University of the Basque Country. PO. Box 644- Bilbao. Spain

*Email: manuel.delasen@ehu.es , santiago.alonso@ehu.es*

** Department of Telecommunications and Systems Engineering

Autonomous University of Barcelona, 08193- Bellaterra, Barcelona, Spain; *Email: Asier.Ibeas@uab.cat*

**Abstract**. This paper is concerned with a SIS (susceptible, infected and susceptible populations) propagation disease model with a nonlinear incidence rate and eventual impulsive (non- necessarily being simultaneous) culling of both populations. The disease transmission does not necessarily take into account the total population as a normalizing effect. In this sense, the considered model is a mixed pseudo-mass. The positivity, stability of both the impulse- free and impulsive under pulse culling variants of the model are investigated in this paper.

**1. Introduction**

Important control problems nowadays related to Life Sciences are the control of ecological models like, for instance, those of population evolution (Beverton-Holt model, Hassell model, Ricker model etc.) via the online adjustment of the species environment carrying capacity, that of the population growth or that of the regulated harvesting quota as well as the disease propagation via vaccination control. In a set of papers, several variants and generalizations of the Beverton-Holt model (like, for instance, standard time–invariant, time-varying parameterized, generalized model or modified generalized model) have been investigated at the levels of stability, cycle- oscillatory behavior, permanence and control through the manipulation of the carrying capacity (see, for instance, [1-5]). The design of related control actions has been proved to be important in those papers at the levels, for instance, of aquaculture exploitation or plague fighting. On the other hand, the literature about epidemic mathematical models is exhaustive in many books and papers including the presence of delays, vaccination strategies and extended impulsive models. A non-exhaustive list of references is given in this manuscript, cf. [6-19] (see also references therein). The sets of epidemic models include the most basic ones, [6-7] as follows:

- SI and SIS- models where not removed- by – immunity population is assumed. In other words, only susceptible and infected populations are assumed. Also, there is no distinction between infected, who do not have external disease symptoms, and infectious who have external symptoms. In SIS models, the infected can recover from illness in contrast to SI models where they cannot.
- SIR models, which include susceptible plus infected plus removed- by –immunity populations without any distinction between infected and infectious.
- SEIR- models where the infected populations is split into two ones (namely, the " infected" which incubate the disease but do not still have any disease symptoms and the " infectious" or " infective" which do have the external disease symptoms).

There are many variants of the above models, for instance, including vaccination of different kinds: constant [8], impulsive [12], discrete – time etc., incorporating point or distributed delays [12-13], periodic or oscillatory behaviors [14] etc. . On the other hand, such models become considerably simpler for the illness transmission among plants [6-7]. The presence of the delays and the positivity of delay-free and delayed dynamic systems have been exhaustively studied in the context of control theory [23-24] including the case of mixed continuous- time and discrete- time hybrid models (see, for instance, [20-22], [25] and references therein). Related issues on discretization and adaptive sampling have also been focused on (see, for instance, [27] and references there in). This paper is concerned with a time-varying SIS (susceptible, infected and susceptible populations) propagation disease model exhibiting a nonlinear incidence rate and impulsive eventual (non-necessarily simultaneous) culling of both populations at a set of time instants separated by an arbitrary distance subject to a minimum threshold in order to guarantee the well-posedness of the model. Pulse culling is an appropriate technique when dealing with animals and plants populations to remove or decrease the infected populations or the exposition of susceptible to infection or both [26]. The nonlinear incidence consists of two time-varying additive terms being respectively proportional to the susceptible and infected populations normalized to the total (which is simultaneously susceptible) population in the absence of disease

## 2. The model and its positivity properties

The susceptible and infected populations are assumed to obey the dynamics:

$$\dot{S}(t) = r(t)\left(1 - \frac{G(t)}{p(t)}\right)S(t) + (\gamma(t) - \beta(t)S(t))I(t) \tag{1.a}$$

$$\dot{I}(t) = [\beta(t)S(t) - d(t) - \gamma(t)]I(t) \tag{1.b}$$

with initial conditions $S(0) = S_0 \geq 0$ and $I(0) = I_0 \geq 0$, where the whole population is $N(t) = S(t) + I(t)$; $\forall t \in \mathbf{R}_{0+}$ with $N(0) = S(0) + I(0) = N_0 = S_0 + I_0 \geq 0$, and

• $PC^{(0)}(\mathbf{R}_{0+}, \mathbf{R}_+)$ and $BPC^{(0)}(\mathbf{R}_{0+}, \mathbf{R}_+)$ are the sets of piecewise continuous and bounded piecewise continuous real functions from $\mathbf{R}_{0+}$ to $\mathbf{R}_+$, respectively. $PC^{(i)}(\mathbf{R}_{0+}, \mathbf{R}_+)$; $i \in \mathbf{Z}_+$ is the set of real functions of class (i-1); i.e. the set of real functions in $C^{(i-1)}(\mathbf{R}_{0+}, \mathbf{R}_+)$ whose i-th derivative exists but it is not necessarily everywhere continuous on its definition domain (thus, $f \in PC^{(i)}(\mathbf{R}_{0+}, \mathbf{R}_+)$ is also in $C^{(i-1)}(\mathbf{R}_{0+}, \mathbf{R}_+)$ but not necessarily in $C^{(i)}(\mathbf{R}_{0+}, \mathbf{R}_+)$). $L^p(\mathbf{R}_{0+}, \mathbf{R})$ is the space of p-integrable real functions of domain $\mathbf{R}_{0+}$.

• $r \in BPC^{(0)}(\mathbf{R}_{0+}, \mathbf{R})$, $d \in BPC^{(0)}(\mathbf{R}_{0+}, \mathbf{R}_{0+})$ and $\gamma \in BPC^{(0)}(\mathbf{R}_{0+}, \mathbf{R}_{0+})$ are the non-identically zero, bounded piecewise continuous intrinsic growth, death and recovery rates of the infected, respectively. The notation for the spaces of functions is as follows:

• $\beta \in BPC^{(0)}(\mathbf{R}_{0+}, \mathbf{R}_+)$ is a periodic infection rate of period $\omega$

• $G \in BPC^{(0)}(\mathbf{R}_{0+}, \mathbf{R}_{0+})$ is the incidence rate which is a piecewise continuous function defined by $G(t) := \delta_1(t)S(t) + \delta_2(t)I(t)$; $\forall t \in \mathbf{R}_{0+}$ where $\delta_i \in BPC^{(0)}(\mathbf{R}_{0+}, \mathbf{R}_{0+})$ with $0 \leq \delta_{mi} \leq \min_{t \in \mathbf{R}_{0+}} \delta_i(t) \leq \max_{t \in \mathbf{R}_{0+}} \delta_i(t) \leq \delta_{Mi} < +\infty$; i=1,2. Note that $\delta_m N(t) \leq G(t) \leq \delta_M N(t)$; $\forall t \in \mathbf{R}_{0+}$ where $\delta_M := \max_{1 \leq i \leq 2}(\delta_{Mi})$ and $\delta_m := \min_{1 \leq i \leq 2}(\delta_{mi}) \leq \delta_M$.

- $p \in BPC^{(0)}(\mathbf{R}_{0+}, \mathbf{R}_+)$ is the carrying capacity function which is defined by $p(t) = K(t) + p_0(t)$; $\forall t \in \mathbf{R}_{0+}$ where $K \in BPC^{(0)}(\mathbf{R}_{0+}, \mathbf{R}_+)$ is subject to $\min_{t \in \mathbf{R}_{0+}} K(t) \geq K_0 + \varepsilon_0 > 0$ for some $K_0, \varepsilon_0 \in \mathbf{R}_+$, $p_0 \in C^{(0)}(\mathbf{R}_{0+}, \mathbf{R})$ is a continuous either oscillating or, in particular, periodic function of period $T_p$ is subject to $\max_{\tau \in [t, t+T_p]} p_0(t) \leq K_0$; $\forall t \in \mathbf{R}_{0+}$. This function can be also used to describe an oscillatory infection-free behavior of the total population. In some results of Sections 4-5, the function $p \equiv p_0 : \mathbf{R}_{0+} \to \mathbf{R}_{0+}$ is allowed to interpret a possible asymptotic extinction of the total population $N(t) = p(t)$ in the absence of infected one. A simple inspection yields that if the functions of parameters of (1) have finite asymptotic limits as $t \to +\infty$, denoted by superscript *, or simply the system is time-invariant then there are two equilibrium points, namely : $S_1^* = I_1^* = 0$ (the total population is extinguished) and the following endemic equilibrium point:

$$S_2^* = \frac{\gamma^* + d^*}{\beta^*} \quad ; \quad I_2^* = \frac{\gamma^*(p^* - G^*)(d^* + \gamma^*)}{d^* \beta^* p^*}$$

The linearized system of the limiting dynamic system about the equilibrium points has the following constant matrix of dynamics:

$$A^* = \begin{bmatrix} r^*\left(1 - 3\frac{\delta_1^* S^*}{p^*}\right) - \left(\beta^* + \frac{2r^* \delta_2^*}{p^*}\right) I^* & \gamma^* - \left(\beta^* + \frac{r^* \delta_2^*}{p^*}\right) S^* \\ \beta^* I^* & -(d^* + \gamma^*) \end{bmatrix}$$

so that $S_1^* = I_1^* = 0 \Rightarrow A_1^* = \begin{bmatrix} r^* & \gamma^* \\ 0 & -(d^* + \gamma^*) \end{bmatrix}$

$$S_2^* = \frac{\gamma^* + d^*}{\beta^*}; \quad I_2^* = \frac{\gamma^*(p^* - \delta_1^* S_2^*)(d^* + \gamma^*)}{d^* \beta^* p^* + r^* \delta_2^*} = \frac{\gamma^*(\beta^* p^* - \delta_1^*(\gamma^* + d^*))(d^* + \gamma^*)}{\beta^*(d^* \beta^* p^* + r^* \delta_2^*)}$$

$$\Rightarrow A_2^* = \begin{bmatrix} r^*\left(1 - 3\frac{\delta_1^*(\gamma^* + d^*)}{p^* \beta^*}\right) - \left(\beta^* + \frac{2r^* \delta_2^*}{p^*}\right) \frac{\gamma^*(\beta^* p^* - \delta_1^*(\gamma^* + d^*))(d^* + \gamma^*)}{\beta^*(d^* \beta^* p^* + r^* \delta_2^*)} & \gamma^* - \left(\beta^* + \frac{r^* \delta_2^*}{p^*}\right) \frac{\gamma^* + d^*}{\beta^*} \\ \frac{\beta^* \gamma^*(\beta^* p^* - \delta_1^*(\gamma^* + d^*))(d^* + \gamma^*)}{\beta^*(d^* \beta^* p^* + r^* \delta_2^*)} & -(d^* + \gamma^*) \end{bmatrix}$$

The subsequent result follows immediately from the discussion of the signs of the real parts or its potential zero values of the eigenvalues of the linearized matrix of dynamics at both equilibrium points:

**Theorem 2.1**. Assume that a limiting system for (2.1) exists with all the parameterizing functions converging to finite limits as time tends to infinity. The following properties hold:

**(i)** The zero equilibrium point $(0,0)^T$ of the limiting system is locally asymptotically stable if and only if $r^* < 0$ and $d^* + \gamma^* > 0$ within an open neighborhood of sufficiently small radius centred at the origin. In this case, the linearized limiting system about the origin is globally exponentially stable.

**(ii)** The zero equilibrium point of the limiting system is locally stable if $r^* = 0$ and $d^* + \gamma^* > 0$ within an open neighborhood of sufficiently small radius centred at the origin. In this case, the linearized system about the origin is globally stable but not globally asymptotically stable. Furthermore, the infected population of the limiting system exponentially extinguishes if the state – trajectory solution enters some neighborhood of the origin of sufficiently small radius. If $r^* = d^* = \gamma^* = 0$ then the limiting system is locally unstable around the zero equilibrium point in an open neighborhood of sufficiently small radius centred at the origin.

**(iii)** The zero equilibrium point of the limiting system is locally unstable if $r^* > 0$ and $d^* + \gamma^* < 0$ within a neighborhood of the origin of sufficiently small radius and the linearized system about the origin is globally unstable.

**(iv)** Let $A_2^* = \left(A_{2ij}^*\right)$ defined by an entry-by-entry matrix notation. Then, the endemic equilibrium point $\left(S_2^*, I_2^*\right)^T$ of the limiting system is locally asymptotically stable within an open neighborhod centred at $\left(S_2^*, I_2^*\right)^T$ of a sufficiently small radius if and only if $A_{211}^* + A_{222}^* < 0$ and $A_{211}^* A_{222}^* > A_{212}^* A_{221}^*$ and the linearized limiting system about $\left(S_2^*, I_2^*\right)^T$ is globally asymptotically stable If $A_{211}^* + A_{222}^* < 0$ and $A_{211}^* A_{222}^* > A_{212}^* A_{221}^*$ or if $A_{211}^* + A_{222}^* > 0$ the endemic equilibrium point of the limiting system about $\left(S_2^*, I_2^*\right)^T$ is locally unstable and that of the linearized limiting system about $\left(S_2^*, I_2^*\right)^T$ is globally unstable. If $A_{211}^* + A_{222}^* \leq 0$ and $A_{211}^* A_{222}^* = A_{212}^* A_{221}^*$ then the endemic equilibrium point of the linearized limiting system is locally stable but not locally asymptotically stable within some neigborhood centred at it of sufficiently small radius. □

*Remark 2.1.* Note that the local stability properties of the limiting system and the global ones of their linearized versions around the equilibrium points are identical if the eigenvalues of the corresponding linearized limiting systems have both positive or negative real values. □

The unique solution of (1.b) is

$$I(t) = I(0) e^{\int_0^t (\beta(\tau) S(\tau) - d(\tau) - \gamma(\tau)) d\tau} \quad ; \forall t \in \mathbf{R}_{0+} \qquad (2.a)$$

which then used in that of (1.a) yields:

$$S(t) = e^{\int_0^t \left[r(\tau)\left(1 - \frac{G(\tau)}{p(\tau)}\right) - \beta(\tau) I(\tau)\right] d\tau} S(0)$$

$$+ I(0) \int_0^t e^{\int_0^\tau \left[ r(\tau') \left( 1 - \frac{G(\tau')}{p(\tau')} \right) - I(0) \beta(\tau') e^{\int_0^{\tau'} (\beta(\tau'') S(\tau'') - d(\tau'') - \gamma(\tau'')) d\tau''} \right] d\tau'}$$

$$\times \gamma(t-\tau) \left( e^{\int_0^{t-\tau} (\beta(\tau') S(\tau') - d(\tau') - \gamma(\tau')) d\tau'} \right) d\tau \; ; \quad \forall t \in \mathbf{R}_{0+} \tag{2.b}$$

$$= S(0) e^{\int_0^t r(\tau) \left( 1 - \frac{G(\tau)}{p(\tau)} \right) d\tau} + \int_0^t e^{\int_0^\tau r(\tau') \left( 1 - \frac{G(\tau')}{p(\tau')} \right) d\tau'} (\gamma(t-\tau) - \beta(t-\tau) S(t-\tau)) I(t-\tau) d\tau \; ;$$

$$\forall t \in \mathbf{R}_{0+} \tag{2.c}$$

The following positivity result is immediate from $S(0) = S_0 \geq 0$ and $I(0) = I_0 \geq 0$, $\gamma(t) \geq 0$ and (2):

**Theorem 2.2** *(Positivity of the solutions of the epidemic model)*. The following properties hold with $S(0) = S_0 \geq 0$ and $I(0) = I_0 \geq 0$, $\gamma(t) \geq 0$:

(i) $S(t) \geq 0$; $I(t) \geq 0$ and $N(t) \geq 0$; $\forall t \in \mathbf{R}_{0+}$.

(ii) $I(0) = 0 \Leftrightarrow (I(t) = 0; \forall t \in \mathbf{R}_{0+})$.

(iii) $N(0) = 0 \Leftrightarrow (S(t) = I(t) = 0; \forall t \in \mathbf{R}_{0+}) \Leftrightarrow (N(t) = 0; \forall t \in \mathbf{R}_{0+})$. □

Theorem 2.2(ii) states that the infected population is zero at any finite time if it is zero at $t=0$. However, this does not imply it can potentially asymptotically extinguish as time tends to infinity for $I(0) > 0$. The fact that the solutions of (1) are nonnegative for all time if their initial conditions are nonnegative imply that the dynamic system (1) is positive. The following stability result follows under Theorem 2.2:

**Corollary 2.3**. The following results hold:

(i) $N(t) = S(t)$ is bounded for $N(0) = S(0)$ being bounded if $g \in L^1(\mathbf{R}_{0+}, \mathbf{R}_{0+})$ where $g(t) = r(t) \left( 1 - \delta_1(t) \frac{N(t)}{p(t)} \right)$. A sufficient condition is $\liminf_{t \to +\infty} \delta_1(t) \geq \liminf_{t \to +\infty} \frac{p(t)}{N(t)}$.

(ii) $N(t) = S(t)$ is bounded for $N(0) = S(0)$ being bounded and converges exponentially to zero as $t \to \infty$ with an upper-bounding function of exponential order of at most $-\varepsilon_{0S} < 0$ for some $\mathbf{R}_+ \ni \varepsilon_S \in (0, \varepsilon_{0S})$ ( so that the population asymptotically extinguishes at an exponential rate ) if

$$\limsup_{t \to +\infty} \left( r(t) \left( 1 - \frac{\delta_1(t) N(t)}{p(t)} \right) + \varepsilon_S t \right) \leq 0. \qquad \text{A sufficient condition is}$$

$$\min \left( \liminf_{t \to \infty} r(t), \liminf_{t \to \infty} \left( \frac{\delta_1(t) N(t)}{p(t)} - 1 \right) \right) \geq \varepsilon_{1S} \qquad \text{for some} \qquad \varepsilon_{1S} \in \mathbf{R}_+. \qquad \text{If}$$

$\lim\limits_{t \to +\infty} \int_0^t r(\tau)\left(1 - \dfrac{\delta_1(\tau)N(\tau)}{p(\tau)}\right)d\tau = -\infty$ then the population asymptotically extinguishes but not necessarily at an exponential rate.

**(iii)** Assume that $S(t_0) = 0$ and $\gamma(t) = 0$; $\forall t \in [t_0, \infty)$ except perhaps on a subset of $[t_0, \infty)$ of zero measure. Then, $S(t) = 0, N(t) = I(t); \forall t \in [t_0, \infty)$. If, in addition, $\gamma(t) = 0$; $\forall t \in [t_0, \infty)$ except perhaps on a subset of $[t_0, \infty)$ of zero measure, then $S(t) = 0, N(t) = I(t) = I(t_0); \forall t \in [t_0, \infty)$. If, furthermore, $\lim\limits_{t \to +\infty} \int_{t_0}^t (d(\tau) + \gamma(\tau))d\tau = +\infty$, then $I(t) \to 0$ asymptotically as $t \to +\infty$. If $\lim\limits_{t \to +\infty} \sup\left(-\int_{t_0}^t (d(\tau) + \gamma(\tau))d\tau + \varepsilon_1 t\right) \leq 0$ then $I(t) \to 0$ asymptotically as $t \to +\infty$ at an exponential rate.

Now, define $x(t) := (S(t), I(t))^T$ and $y(t) = N(t)$, so that the differential system (1) becomes the equivalent dynamic system of state vector function $x: \boldsymbol{R}_{0+} \to \boldsymbol{R}^2$ and output function $y: \boldsymbol{R}_{0+} \to \boldsymbol{R}$

$$\dot{x}(t) = A(t, x(\boldsymbol{\tau}_t))x(t, x(\boldsymbol{\tau}_t)) = A_0(t, x\boldsymbol{\tau}_t)x(t, x(\boldsymbol{\tau}_t)) \; ; \; y(t, x(\boldsymbol{\tau}_t)) = c^T x(t, x(\boldsymbol{\tau}_t)) \qquad (3)$$

with initial conditions $x(0) = x_0 = (S(0), I(0))^T = (S_0, I_0)^T$ with $\boldsymbol{\tau}_t$ being the time interval $[0, t)$ where

$$A(t, x(\boldsymbol{\tau}_t)) := \begin{bmatrix} r(t)\left(1 - \dfrac{G(t)}{p(t)}\right) - \beta(t)I(t) & \gamma(t) \\ 0 & \beta(t)S(t) - d(t) - \gamma(t) \end{bmatrix} ; \; c := (1, 1)^T \qquad (4.a)$$

$$A_0(t, x(\boldsymbol{\tau}_t)) := \begin{bmatrix} 0 & r(t)\left(1 - \dfrac{G(t)}{p(t)}\right) - \beta(t)S(t) + \gamma(t) \\ 0 & \beta(t)S(t) - d(t) - \gamma(t) \end{bmatrix} \qquad (4.b)$$

whose solution is for initial condition $x(0) = x_0$ after direct combination of $S(t)$ and $I(t)$:

$$x(t) = \Psi(t, x(\boldsymbol{\tau}_t))x(0) = \begin{bmatrix} \Psi_{11}(t, x(\boldsymbol{\tau}_t)) & \Psi_{12}(t, x(\boldsymbol{\tau}_t)) \\ 0 & \Psi_{22}(t, x(\boldsymbol{\tau}_t)) \end{bmatrix}x(0) \qquad (5.a)$$

$$= \Psi_0(t, x(\boldsymbol{\tau}_t))x(0) = \begin{bmatrix} 0 & \Psi_{012}(t, x(\boldsymbol{\tau}_t)) \\ 0 & \Psi_{022}(t, x(\boldsymbol{\tau}_t)) \end{bmatrix}x(0) \qquad (5.b)$$

; $\forall t \in \boldsymbol{R}_{0+}$, with $A: \boldsymbol{R}_{0+} \times \boldsymbol{R}^n \to \boldsymbol{R}^n$ and $A_0: \boldsymbol{R}_{0+} \times \boldsymbol{R}^n \to \boldsymbol{R}^n$ being two matrix functions of dynamics describing equivalently the dynamic system, and $\Psi \in L(\boldsymbol{R}^2, \boldsymbol{R}^2)$ and $\Psi_0 \in L(\boldsymbol{R}^2, \boldsymbol{R}^2)$ are the nonlinear operators whose representations are the matrix functions $\Psi(t, x(\boldsymbol{\tau}_t))$ and $\Psi_0(t, x(\boldsymbol{\tau}_t))$, respectively. Note that the nonlinear operators $\Psi: \boldsymbol{R}^2 \to \boldsymbol{R}^2$ and $\Psi_0: \boldsymbol{R}^2 \to \boldsymbol{R}^2$ are not necessarily bounded (and then non-necessarily continuous) in the normed Euclidean space $\boldsymbol{R}^2$ since the entries of the associated fundamental matrix function of (1), $\Psi \in PC^{(1)}(\boldsymbol{R}_{0+}, \boldsymbol{R}^2)$ can eventually be non-continuous and unbounded depending on its parameterizing functions, then the dynamic system (3) being unstable.

The fundamental matrix functions $\Psi \in PC^{(1)}(\mathbf{R}_{0+}, \mathbf{R}^2)$ and $\Psi_0 \in PC^{(1)}(\mathbf{R}_{0+}, \mathbf{R}^2)$ may also be viewed as evolution operators of the dynamic system as reflected by (3), with:

$$\Psi_{11}(t, x(\tau_t)) := \int_0^t \left[ r(\tau)\left(1 - \frac{G(\tau)}{p(\tau)}\right) - I(0)\beta(\tau) e^{\int_0^\tau (\beta(\tau')S(\tau') - d(\tau') - \gamma(\tau'))d\tau'} \right] d\tau$$

$$\Psi_{12}(t, x(\tau_t)) := \int_0^t e^{\int_0^\tau \left[ r(\tau')\left(1 - \frac{G(\tau')}{p(\tau')}\right) - I(0)\beta(\tau') e^{\int_0^{\tau'} (\beta(\tau'')S(\tau'') - d(\tau'') - \gamma(\tau''))d\tau''} \right] d\tau'} d\tau$$

$$\Psi_{22}(t, x(\tau_t)) := e^{\int_0^t (\beta(\tau)S(\tau) - d(\tau) - \gamma(\tau))d\tau} \tag{6.a}$$

$$\Psi_{012}(t, x(\tau_t)) := e^{\int_0^t \left[ r(\tau)\left(1 - \frac{G(\tau)}{p(\tau)}\right) - \beta(\tau)S(\tau) + \gamma(\tau) \right] d\tau} \quad ; \quad \Psi_{022}(t, x(\tau_t)) = \Psi_{22}(t, x(\tau_t)) \tag{6.b}$$

; $\forall t \in \mathbf{R}_{0+}$. The following result is a direct consequence of Theorem 2.2.

**Corollary 2.4** ( *Positivity of the solutions of the dynamic system (3) –(4)*). The dynamic system (3) –(4) is a positive dynamic system in the sense that both state components are nonnegative for all time if both components of the initial state are nonnegative for all time. □

Assume that the recovery rate function of the infected, $\gamma$, is allowed to be negative on some (perhaps non-connected) time intervals of finite measure. The interpretation is that the susceptible rate decreases and the infected rate increases within such intervals what may be interpreted by an increase of the disease pathogens due to atmospheric conditions, social behaviours of susceptible and infected, etc. Consider the disjoint union of nonnegative real intervals $\mathbf{R}_{0+} := R_{\gamma+} \cup R_{\gamma-} \cup R_{\gamma 0}$ with $R_{\gamma-}$, $R_{\gamma 0}$ being empty or non-empty where.

$$R_{\gamma+} := \{z \in \mathbf{R}_{0+} : \gamma(z) > 0\}; \quad R_{\gamma-} := \{z \in \mathbf{R}_{0+} : \gamma(z) < 0\}; \quad R_{\gamma 0} := \{z \in \mathbf{R}_{0+} : \gamma(z) = 0\}$$

and also

$$R_{t\gamma+} := \{z(<t) \in \mathbf{R}_{0+} : \gamma(t) > 0\}; \quad R_{t\gamma-} := \{z(<t) \in \mathbf{R}_{0+} : \gamma(t) < 0\}; \quad R_{t\gamma 0} := \{z(<t) \in \mathbf{R}_{0+} : \gamma(t) = 0\}; \quad \forall t \in \mathbf{R}_+$$

such that $\mathbf{R}_{0+} \supset [0, t) = R_{t\gamma+} \cup R_{t\gamma-} \cup R_{t\gamma 0}$. Eqn. (2.b) becomes:

$$S(t) = e^{\int_0^t \left[ r(\tau)\left(1 - \frac{G(\tau)}{p(\tau)}\right) - \beta(\tau)I(\tau) \right] d\tau} S(0)$$

$$+ I(0) \int_{R_{t\gamma+}} e^{\int_0^\tau \left[ r(\tau')\left(1 - \frac{G(\tau')}{p(\tau')}\right) - I(0)\beta(\tau') e^{\int_0^{\tau'} (\beta(\tau'')S(\tau'') - d(\tau'') - \gamma(\tau''))d\tau''} \right] d\tau'} d\tau'$$

$$\times \gamma(t-\tau) \left( e^{\int_0^{t-\tau} (\beta(\tau') S(\tau') - d(\tau') - \gamma(\tau')) d\tau'} \right) d\tau$$

$$- I(0) \int_{R_{t\gamma+}} e^{\int_0^\tau \left[ r(\tau') \left( 1 - \frac{G(\tau')}{p(\tau')} \right) - I(0) \beta(\tau') e^{\int_0^{\tau'} (\beta(\tau'') S(\tau'') - d(\tau'') - \gamma(\tau'')) d\tau''} \right] d\tau'}$$

$$\times |\gamma(t-\tau)| \left( e^{\int_0^{t-\tau} (\beta(\tau') S(\tau') - d(\tau') - \gamma(\tau')) d\tau'} \right) d\tau \; ; \; \forall t \in \mathbf{R}_{0+} \qquad (7)$$

The following result is immediate from Theorem 2.1 (i), Corollary 2.4, (2.a) and (7):

**Corollary 2.5** *(Positivity of the solutions of the epidemic model and those of the dynamic system (3)-(4) under eventual negative values of the recovery rate of the infected ).* Theorem 2.2 (i) and Corollary 2.3 still hold if $I(0) > 0$ and any of the two conditions below hold:

(i) $\int_{R_{t\gamma+}} e^{\int_0^\tau \left[ r(\tau') \left( 1 - \frac{G(\tau')}{p(\tau')} \right) - I(0) \beta(\tau') e^{\int_0^{\tau'} (\beta(\tau'') S(\tau'') - d(\tau'') - \gamma(\tau'')) d\tau''} \right] d\tau'}$

$$\times |\gamma(t-\tau)| \left( e^{\int_0^{t-\tau} (\beta(\tau') S(\tau') - d(\tau') - \gamma(\tau')) d\tau'} \right) d\tau$$

$$\leq \left\{ e^{\int_0^t \left[ r(\tau) \left( 1 - \frac{G(\tau)}{p(\tau)} \right) - \beta(\tau) I(\tau) \right] d\tau} \left( S(0)/I(0) \right) \right.$$

$$+ \int_{R_{t\gamma+}} e^{\int_0^\tau \left[ r(\tau') \left( 1 - \frac{G(\tau')}{p(\tau')} \right) - I(0) \beta(\tau') e^{\int_0^{\tau'} (\beta(\tau'') S(\tau'') - d(\tau'') - \gamma(\tau'')) d\tau''} \right] d\tau'}$$

$$\left. \times \gamma(t-\tau) \left( e^{\int_0^{t-\tau} (\beta(\tau') S(\tau') - d(\tau') - \gamma(\tau')) d\tau'} \right) d\tau \right\} \; ; \; \forall t \in \mathbf{R}_{0+}$$

(ii) $\sup_{\tau \in R_{t\gamma-}} (|\gamma(\tau)|) \leq \left\{ \int_{R_{t\gamma+}} e^{\int_0^\tau \left[ r(\tau') \left( 1 - \frac{G(\tau')}{p(\tau')} \right) - I(0) \beta(\tau') e^{\int_0^{\tau'} (\beta(\tau'') S(\tau'') - d(\tau'') - \gamma(\tau'')) d\tau''} \right] d\tau'} \right.$

$$\left. \times \left( e^{\int_0^{t-\tau} (\beta(\tau') S(\tau') - d(\tau') - \gamma(\tau')) d\tau'} \right) d\tau \right\}^{-1}$$

$$\times \left( e^{\int_0^t \left[ r(\tau)\left(1-\frac{G(\tau)}{p(\tau)}\right)-\beta(\tau)I(\tau)\right]d\tau} \right) (S(0)/I(0))$$

$$+ \int_{R_{t\gamma+}} e^{\int_0^\tau \left[ r(\tau')\left(1-\frac{G(\tau')}{p(\tau')}\right)-I(0)\,\beta(\tau')\,e^{\int_0^{\tau'}(\beta(\tau'')S(\tau'')-d(\tau'')-\gamma(\tau''))d\tau''}\right]d\tau'}$$

$$\times \gamma(t-\tau)\left( e^{\int_0^{t-\tau}(\beta(\tau')S(\tau')-d(\tau')-\gamma(\tau'))d\tau'} \right) d\tau \quad ; \forall t \in R_{0+} \quad \square$$

## 3. Stability and comparison to the infection- free case

The following technical assumptions are introduced, [26]:

*Assumption 3.1.* In the absence of disease, the population grows in logistic fashion with positive time-varying carrying capacity p(t) and an intrinsic growth rate function r(t). $\square$

*Assumption 3. 2.* In the presence of disease, only the susceptible population contributes to its carrying capacity. $\square$

Summing- up both sides of both Eqs. 1 leads to the dynamics of the whole population $N(t)=S(t)+I(t)$ as follows by taking into account that $G(t):=\delta_1(t)S(t)+\delta_2(t)I(t)$; $\forall t \in R_{0+}$:

$$r(t)\left(1-\frac{\delta_M N(t)}{p(t)}\right)S(t)-d(t)I(t)$$

$$\leq \dot{N}(t)=r(t)\left(1-\frac{G(t)}{p(t)}\right)S(t)-d(t)I(t)\leq r(t)\left(1-\frac{\delta_m N(t)}{p(t)}\right)S(t)-d(t)I(t)$$

(7.a)

$$=r(t)\left(1-\frac{G(t)}{p(t)}\right)N(t)-\left(r(t)\left(1-\frac{G(t)}{p(t)}\right)+d(t)\right)I(t) \tag{7.b}$$

$$=r(t)\left(1-\frac{\delta_1(t)S(t)}{p(t)}\right)S(t)-\left(r(t)\frac{\delta_2(t)S(t)}{p(t)}+d(t)\right)I(t) \tag{7.c}$$

$$=r(t)\left(1-\frac{\delta_1(t)(N(t)-I(t))}{p(t)}\right)N(t)+\left(r(t)\left(\frac{(\delta_1(t)-\delta_2(t))(N(t)-I(t))}{p(t)}-1\right)-d(t)\right)I(t) \tag{7.d}$$

$$=r(t)N(t)\left(1-\frac{\delta_1(t)N(t)}{p(t)}\right)+\frac{(2\delta_1(t)-\delta_2(t))r(t)N(t)I(t)}{p(t)}-(r(t)+d(t))I(t))$$

$$-\frac{r(t)(\delta_1(t)-\delta_2(t))I^2(t)}{p(t)} \tag{7.e}$$

$$= r(t)N(t)\left(1 - \frac{\delta_1(t)N(t)}{p(t)}\right) + \frac{r(t)}{p(t)}\left[(N(t)+S(t))\delta_1(t) - S(t)\delta_2(t) - (r(t)+d(t))\frac{p(t)}{r(t)}\right]I(t)$$

Since

$$\frac{\delta_2(t)S(t) + p(t)(1+d(t)/r(t))}{N(t)+S(t)} \geq \delta_1(t) \geq \frac{p(t)}{N(t)} \quad \text{then } \dot{N}(t) \leq 0 \text{ and } N(t) \leq N(0); \forall t \in \mathbf{R}_{0+} \text{ with}$$

$N: \mathbf{R}_{0+} \to \mathbf{R}_{0+}$ being monotone non-increasing everywhere in $\mathbf{R}_{0+}$ which has the subsequent solution after using the solution (2.a) of (1.b):

$$N(t) = N(0)e^{\int_0^t r(\tau)\left(1-\frac{G(\tau)}{p(\tau)}\right)d\tau} - I(0)\int_0^t e^{\int_0^\tau r(\tau')\left(1-\frac{G(\tau')}{p(\tau')}\right)d\tau'} \left(r(t-\tau)\left(1-\frac{G(t-\tau)}{p(t-\tau)}\right) + d(t-\tau)\right)$$

$$\times \left(\int_0^{t-\tau} (\beta(\tau')S(\tau') - d(\tau') - \gamma(\tau'))d\tau'\right)d\tau \tag{8.a}$$

$$= N(0)e^{\int_0^t r(\tau)\left(1-\frac{\delta_1(\tau)S(\tau)}{p(\tau)}\right)d\tau}$$

$$+ \int_0^t e^{\int_0^\tau r(\tau')\left(1-\frac{\delta_1(\tau')S(\tau')}{p(\tau')}\right)d\tau'} \left(r(t-\tau)\left(\frac{(\delta_1(t-\tau)+\delta_2(t-\tau))S(t-\tau)}{p(t-\tau)} - 1\right) - d(t-\tau)\right)I(t-\tau)d\tau$$

$$\tag{8.b}$$

$$= N(0)e^{\int_0^t r(\tau)\left(1-\frac{\delta_1(\tau)S(\tau)}{p(\tau)}\right)d\tau}$$

$$+ I(0)\int_0^t e^{\int_0^\tau r(\tau')\left(1-\frac{\delta_1(\tau')S(\tau')}{p(\tau')}\right)d\tau'} \left(r(t-\tau)\left(\frac{(\delta_1(t-\tau)+\delta_2(t-\tau))S(t-\tau)}{p(t-\tau)} - 1\right) - d(t-\tau)\right)$$

$$\times e^{\int_0^{t-\tau}(\beta(\tau')S(\tau') - d(\tau') - \gamma(\tau'))d\tau'} d\tau \; ; \; \forall t \in \mathbf{R}_{0+} \tag{8.c}$$

It is interesting to investigate the boundedness of the solutions of the system (1) when it is positive, its intrinsic growth and death rates are strictly positive for all time and $\delta_1 \in BPC^{(0)}(\mathbf{R}_{0+}, \mathbf{R}_+)$ which are real constraints in practical situations. This is addressed in the subsequent result.

**Theorem 3.1**. Assume that the system (1) is positive with nonnegative bounded initial conditions so that it satisfies Theorem 2.2. Assume also that $\delta_{m1} > 0$, $d \in BPC^{(0)}(\mathbf{R}_{0+}, \mathbf{R}_+)$ and $r \in BPC^{(0)}(\mathbf{R}_{0+}, \mathbf{R}_+)$. Then, the susceptible, infected and total populations are bounded for all time for any given set of bounded nonnegative initial conditions $S(0)$ and $I(0)$. □

It is interesting to see that even if $G(t) = 0$ and $r(t) > 0$; $\forall t \in \mathbf{R}_{0+}$ the total, susceptible and infected populations are not unbounded as addressed in the subsequent result:

**Theorem 3.2**. Assume that $\liminf_{t\to+\infty} \min(d(t), r(t)) > 0$. Then, any nonnegative solution of (1) has the total, susceptible and infected populations bounded for all time under bounded initial conditions if $\liminf_{t\to+\infty} d(t)/r(t) > 0$ is sufficiently large. The property also holds if $\liminf_{t\to+\infty} d(t) > 0$ and $\exists \lim_{t\to+\infty} r(t) = 0$. □

Some further boundedness results for the populations of the model (1) follow below:

**Theorem 3.3**. The following properties hold:

(i) Assume that $\lim_{t\to+\infty} \int_0^t r(\tau)\left(1 - \frac{G(\tau)}{p(\tau)}\right) d\tau = -\infty$ and $\liminf_{t\to+\infty}\left(d(\tau) + r(\tau)\left(1 - \frac{G(\tau)}{p(\tau)}\right)\right) \geq 0$. A sufficient condition for those constraints to hold is $G(t) \leq K_N(t)N(t)$; $\forall t \in R_{0+}$, guaranteed if $K_N(t) = \max(\delta_1(t), \delta_2(t))$; $\forall t \in R_{0+}$, with the subsequent subset of $R_{0+}$ being of zero measure:

$Z_G := \{t \in R_{0+} : (p(t)) \leq G(t) \vee r(t) \leq 0\}$

Then, $N(t) \to 0$, $S(t) \to 0$ and $I(t) \to 0$ asymptotically as $t \to +\infty$, and

$$\lim_{t\to+\infty} \int_0^t e^{\int_0^\tau r(\tau')\left(1 - \frac{G(\tau')}{p(\tau')}\right)d\tau'} \left[r(t-\tau)\left(1 - \frac{G(t-\tau)}{p(t-\tau)}\right) + d(t-\tau)\right]\left(e^{\int_0^{t-\tau}(\beta(\tau')S(\tau') - d(\tau') - \gamma(\tau'))d\tau'}\right) d\tau = 0$$

if $I(0) > 0$. (9)

(ii) Assume that $\limsup_{t\to+\infty} S(t) < +\infty$ and there exist $\varepsilon_I \in R_+$, $\varepsilon_d, \varepsilon_\gamma \in R_{0+}$ such that

$R_+ \ni (\varepsilon_d + \varepsilon_\gamma) \geq \varepsilon_I$, $\liminf_{t\to+\infty} d(t) \geq \varepsilon_d$, $\liminf_{t\to+\infty} \gamma(t) \geq \varepsilon_\gamma$ and $\limsup_{t\to+\infty} \beta(t) \leq \frac{\varepsilon_d + \varepsilon_\gamma - \varepsilon_I}{\limsup_{t\to+\infty} S(t)}$. Then,

$I : R_{0+} \to R_{0+}$ is bounded and globally exponentially stable to the origin with an upper-bounding function of exponential order of at most $-\varepsilon_I < 0$.

(ii.1) Assume that $\delta_2(t) < \min\left(1, \frac{p(t)d(t)}{r(t)S(t)}\right)$; $\forall t(\geq \bar{t}) \in R_{0+}$ for some finite $\bar{t} \in R_{0+}$ (a sufficient condition being $\delta_2(t) < \min\left(1, \frac{p(t)d(t)}{r(t)N(t)}\right)$; $\forall t(\geq \bar{t}) \in R_{0+}$). Then, $N : R_{0+} \to R_{0+}$ is bounded.

(ii.2) Assume that $\delta_2(t_i) < \min\left(1, \frac{p(t_i)d(t_i)}{r(t_i)S(t_i)}\right)$ at the real sequence $ST$ of elements

$t_i := \{\min \tau : R_{0+} \ni \tau \geq \max(\bar{t}, t_{i-1}) : N(t_i) \geq \overline{N}\}$; $i \in IST \subset Z_+$ ( $IST$ being the indexing set of ST) for any prefixed triple $\bar{t}, \overline{N}, \overline{N}_\varepsilon(<\overline{N}) \in R_+$, such that $N(t) = \overline{N} + \overline{K}_{\overline{N}} \overline{K}_{\overline{N}t_i}$ for some existing bound real sequence $\overline{K}_{\overline{N}t_i} = \overline{K}_{\overline{N}t_i}(\bar{t}, \overline{N}, \overline{N}_\varepsilon, t_i) \in R_{0+}$. Assume also that

$$\delta_2(t_i+\alpha) < min\left(1, \frac{p(t_i+\alpha)d(t_i+\alpha)}{r(t_i+\alpha)S(t_i+\alpha)}\right); \quad \forall \alpha \in [0,\alpha_{0i}] \subset \mathbf{R}_{0+}, \quad \forall i \in IST, \text{ for some existing}$$

bounded sequence $\alpha_{0i} \in \mathbf{R}_+$, such that $\overline{N} - \overline{N}_\varepsilon \leq N(t_i + \alpha) \leq \overline{N} + \overline{K}_{\overline{N}t_i}$; $\forall t_i \in ST$. Thus,

$$\overline{N} - \overline{N}_\varepsilon \leq N(t) \leq \overline{N} + \overline{K}_{\overline{N}t} < \infty; \quad \forall t \geq \overline{t} \text{ and } \sup_{t \in \mathbf{R}_{0+}} N(t) < \infty. \qquad \square$$

**Assertion 3.4**. Assume that $\limsup\limits_{t \to +\infty}\left(\frac{d(t)}{r(t)} - \frac{S(t)}{I(t)}\right) \leq 0$, equivalently $\limsup\limits_{t \to +\infty}\left(1 + \frac{d(t)}{r(t)} - \frac{N(t)}{I(t)}\right) \geq 0$, and $\liminf\limits_{t \to +\infty} \frac{d(t)}{r(t)} \geq 0$. Then $\limsup\limits_{t \to +\infty} N(t) < +\infty$ and $N(t) \leq \overline{N} < +\infty$; $\forall t \in \mathbf{R}_{0+}$, with $\overline{N}$ being dependent on initial conditions, for any given nonnegative initial conditions.

**Proof:** Assume that $\limsup\limits_{t \to +\infty}\left(\frac{d(t)}{r(t)} - \frac{S(t)}{I(t)}\right) \leq 0$, equivalently $\limsup\limits_{t \to +\infty}\left(1 + \frac{d(t)}{r(t)} - \frac{N(t)}{I(t)}\right) \leq 0$, and

$\liminf\limits_{t \to +\infty} \frac{d(t)}{r(t)} \geq 0$ and proceed by contradiction by also assuming $\lim\limits_{t \to +\infty} N(t) = +\infty$. It turns out that it exits some sufficiently large finite $M \in \mathbf{R}_+$ such that $\dot{N}(t) \leq 0$ if $N(t) \geq M$ since

$$\delta_M(t) \geq \delta_m(t) \geq max\left(0, \frac{p(t)}{N(t)}\left(1 - \frac{d(t)I(t)}{r(t)S(t)}\right)\right) \Rightarrow \dot{N}(t) \leq 0; \quad \forall t \in \mathbf{R}_{0+}$$

so that the $N: \mathbf{R}_{0+} \to \mathbf{R}_{0+}$ is non-increasing on subset $\{t \in \mathbf{R}_{0+} : N(t) \geq M\}$ of $\mathbf{R}_{0+}$ is of a finite measure depending on initial conditions (in particular, it is empty if $N(0) \leq M$). This contradicts $N: \mathbf{R}_{0+} \to \mathbf{R}_{0+}$ and then $N: \mathbf{R}_{0+} \to \mathbf{R}_{0+}$ is bounded for any given set of nonnegative bounded initial conditions. $\qquad \square$

**Theorem 3.5**. The following properties hold:

**(i)** If $\int_0^t r(\tau)\left(1 - \frac{\delta_m N(\tau)}{p(t)}\right)d\tau < +\infty$; $\forall t \in \mathbf{R}_{0+}$ then $N: \mathbf{R}_{0+} \to \mathbf{R}_{0+}$, its time derivative and the partial populations of susceptible and infected are all bounded for any bounded initial conditions. A sufficient condition is $\int_0^t r(\tau)\,d\tau < +\infty$; $\forall t \in \mathbf{R}_{0+}$. Another weaker sufficient condition is $\int_0^t r(\tau)\vartheta(\tau)\,d\tau < +\infty$; $\forall t \in \mathbf{R}_{0+}$ where $\vartheta : \mathbf{R}_{0+} \to \{0,1\}$ is a binary indicator function defined as $\vartheta(\tau) = 1$ if $N(\tau) < p(\tau)/\delta_m$ and $\vartheta(\tau) = 0$, otherwise.

**(ii)** If $d: \mathbf{R}_{0+} \to \mathbf{R}_{0+}$ then the sufficient conditions of (i) guarantee that $\int_0^t d(\tau)I(\tau)d\tau < +\infty$; $\forall t \in \mathbf{R}_{0+}$.

**(iii)** If $\dfrac{d(t)}{r(t)} \geq \dfrac{G(t)-p(t)}{p(t)}$; $\forall t \in \mathbf{R}_{0+}$ then the sufficient conditions of (i) guarantee that

$$\int_0^t \left[ r(\tau)\left(1-\dfrac{G(\tau)}{p(\tau)}\right) + d(\tau) \right] I(\tau)d\tau < +\infty; \ \forall t \in \mathbf{R}_{0+}.$$

**(iv)** If $d: \mathbf{R}_{0+} \to \mathbf{R}$ then the sufficient conditions of (i) guarantee that $\int_0^t d(\tau)\vartheta_d(\tau)I(\tau)d\tau < +\infty$; $\forall t \in \mathbf{R}_{0+}$, where $\vartheta_d: \mathbf{R}_{0+} \to \{0,1\}$ is a binary indicator function defined as $\vartheta_d(\tau) = 1$ if $d(\tau) > 0$ and $\vartheta_d(\tau) = 0$, otherwise.

**(v)** If $\dfrac{d(t)}{r(t)} \geq \dfrac{G(t)-p(t)}{p(t)}$; $\forall t \in \mathbf{R}_{0+}$ then the sufficient conditions of (i) guarantee that

$$\int_0^t \left[ r(\tau)\left(1-\dfrac{G(\tau)}{p(\tau)}\right) + d(\tau) \right] I(\tau)d\tau < +\infty, \ \forall t \in \mathbf{R}_{0+} \text{ where } \vartheta_{dr}: \mathbf{R}_{0+} \to \{0,1\} \text{ is a binary}$$

indicator function defined as $\vartheta_d(\tau) = 1$ if $\dfrac{d(t)}{r(t)} \geq \dfrac{G(t)-p(t)}{p(t)}$ and $\vartheta_{dr}(\tau) = 0$, otherwise. □

Note that:

$$N(t)\dot{N}(t) = r(t)\left(1-\dfrac{G(t)}{p(t)}\right)S(t)N(t) - d(t)I(t)N(t)$$

$$= r(t)\left(1-\dfrac{G(t)}{p(t)}\right)N^2(t) - \left[r(t)\left(1-\dfrac{G(t)}{p(t)}\right) + d(t)\right]I(t)N(t)$$

This result will be useful to prove the following result:

**Theorem 3.6.** If $0 \leq \dfrac{d(t)}{r(t)} < \dfrac{S(t)}{I(t)}$; $\forall t \in \mathbf{R}_{0+}$ then $S^* = I^* = 0$ is a locally unstable equilibrium point. If $\dfrac{d(t)}{r(t)} > \dfrac{S(t)}{I(t)}$; $\forall t \in \mathbf{R}_{0+}$ then $S^* = I^* = 0$ is a locally asymptotically stable equilibrium point. □

It follows from Theorem 2.2, Assertions 3.4 and Theorem 3.6 that the global Lyapunov stability property with ultimate boundedness might be compatible with the existence of a locally unstable zero equilibrium point as follows:

**Theorem 3.7.** Assume that $\liminf\limits_{t \to +\infty} \dfrac{d(t)}{r(t)} \geq 0$ and $\limsup\limits_{t \to +\infty} \left(\dfrac{d(t)}{r(t)} - \dfrac{S(t)}{I(t)}\right) \leq 0$. Then, the SIS model is globally Lyapunov stable for any given bounded nonnegative initial conditions and exhibits the ultimate boundedness property. Furthermore, if $\dfrac{d(t)}{r(t)} \leq \dfrac{S(t)}{I(t)}$; $\forall t \in \mathbf{R}_{0+}$ then $S^* = I^* = 0$ is locally unstable around the equilibrium point. □

The following result is a direct consequence of (2.a), (2.b), where its two right hand-side terms are nonnegative for nonnegative initial conditions, and Theorem 2.2 since the uniform boundedness of the total population implies that of the susceptible and infected ones for all nonnegative system solutions under Theorem 2.2.

**Theorem 3.8.** If $\limsup_{t \to +\infty} \left( \frac{d(t)}{r(t)} - \frac{S(t)}{I(t)} \right) \leq 0$ and $\liminf_{t \to +\infty} \frac{d(t)}{r(t)} \geq 0$ then

**(i)** $-\infty \leq \int_0^\infty (\beta(\tau)S(\tau) - d(\tau) - \gamma(\tau)) d\tau < +\infty$ ; $\limsup_{t \to +\infty} (\beta(t)S(t) - d(t) - \gamma(t)) \leq 0$

and either it exists a finite $T_1 \in \mathbf{R}_{0+}$ such that $\beta(t) = \frac{d(t) + \gamma(t)}{S(t)}$ ; $\forall t \in [T_1, \infty) \setminus S1$ for some empty or nonempty subset $S1$ of zero measure of $[T_1, \infty)$ or $\limsup_{t \to +\infty} (\beta(t)S(t) - d(t) - \gamma(t)) < 0$.

**(ii)** $-\infty \leq \int_0^\infty \left[ r(\tau)\left(1 - \frac{G(\tau)}{p(\tau)}\right) - \beta(\tau)I(\tau) \right] d\tau < +\infty$ ; $\limsup_{t \to +\infty} \left[ r(t)\left(1 - \frac{G(t)}{p(t)}\right) - \beta(t)I(t) \right] \leq 0$

and either it exists a finite $T_2 \in \mathbf{R}_{0+}$ such that $\beta(t) = \frac{r(t)}{I(t)}\left(1 - \frac{G(t)}{p(t)}\right)$ ; $\forall t \in [T_2, \infty) \setminus S2$ for some empty or nonempty subset $S2$ of zero measure of $[T_2, \infty)$ or $\limsup_{t \to +\infty} \left[ r(t)\left(1 - \frac{G(t)}{p(t)}\right) - \beta(t)I(t) \right] < 0$.

**(iii)** If $\liminf_{t \to +\infty} (\beta(t)S(t) - d(t) - \gamma(t)) > 0$ then both the infected and the total populations are unbounded. □

## 4. About the potential absence of susceptible population and the potential absence of infection

It is now proven that, except in trivial cases, the susceptible population cannot be identically zero while the infected one can be identically zero in the SIS- model (1). The second case represents the absence of infection in the SIS–epidemic model (1). This situation is a key tool in fixing values for a periodic, or in particular oscillatory, function $p: \mathbf{R}_{0+} \to \mathbf{R}_+$ in the epidemic SIS- model (1) so that the whole population $N(t) = S(t)$; $\forall t \in \mathbf{R}_{0+}$ evolves through time in a periodic way in the absence of disease.

**Assertion 4.1.** Assume that $\gamma: \mathbf{R}_{0+} \to \mathbf{R}_{0+}$ is nonzero within some subset of nonzero finite measure of its definition domain. Then, the SIS- model (1) with identically zero susceptible population is unfeasible except for its trivial solution. □

**Theorem 4.2.** Assume that the infected population is identically zero, that $\delta_1(t)$ is non identical to $p(t)/N(t)$ except perhaps within some set of finite measure and $r \in BPC^{(0)}(\mathbf{R}_{0+}, \mathbf{R}_+)$ except perhaps on some set of zero measure. Then, the total population equalizing the susceptible one is an oscillatory bounded function.

**Proof**: If $I \equiv 0$ then, one gets from (1):

$$\dot{N}(t) = \dot{S}(t) = r(t)\left(1 - \frac{\delta_1(t)N(t)}{p(t)}\right)N(t) \; ; \; N(0) = N_0 > 0'$$

so that by decomposing $[0, t) = R_{\delta t+} \cup R_{\delta t-} \cup R'_{\delta t}$ ; $\forall t \in R_{0+}$, where

$$R_{\delta t+} := \left\{ R_{0+} \ni \tau (\le t) : \left(\frac{\delta_1(\tau)N(\tau)}{p(\tau)}\right) < 1 \land r(t) \in R_+ \right\}$$

$$R_{\delta t-} := \left\{ R_{0+} \ni \tau (\le t) : \left(\frac{\delta_1(\tau)N(\tau)}{p(\tau)}\right) \ge 1 \land r(t) \in R_+ \right\}$$

$$R'_{\delta t} := \left\{ R_{0+} \ni \tau (\le t) : r(t) \notin R_+ \right\}$$

; $\forall t \in R_{0+}$ are disjoint and (in general) non-connected real intervals with $R'_{\delta t}$ being of finite Lebesgue measure. Then, the solution of the total population is given by:

$$0 \le N(t) = N(0) + \int_{R_{\delta t+}} r(\tau)\left(1 - \frac{\delta_1(\tau)N(\tau)}{p(\tau)}\right)N(\tau)d\tau - \int_{R_{\delta t-}} r(\tau)\left|1 - \frac{\delta_1(\tau)N(\tau)}{p(\tau)}\right|N(\tau)d\tau$$

$$+ \int_{R'_{\delta t}} r(\tau)\left(1 - \frac{\delta_1(\tau)N(\tau)}{p(\tau)}\right)N(\tau)d\tau \; ; \; \forall t \in R_{0+}$$

Proceed by contradiction. If $N : R_{0+} \to R_{0+}$ is unbounded then there is a strictly increasing real sequence $\{t_k \in R_{0+}\}_{k \in Z_{0+}}$ such that $N(t_k) \to +\infty$ as $R_{0+} \ni t_k \to +\infty$ (so as $Z_{0+} \ni k \to +\infty$) and $\lim_{k \to +\infty} \mu(R_{\delta t_k +}) = +\infty$ where $\mu(R_{\delta t_k +})$ is the Lebesgue measure of the set $R_{\delta t_k +} = \left\{ R_{0+} \ni \tau(\le t_k) : \left(\frac{\delta_1(\tau)N(\tau)}{p(\tau)}\right) < 1 \land r(t) \in R_+ \right\}$. But then $0 < \frac{\delta_{1m} N(t_k)}{p(t_k)} (\to +\infty) < 1$ as $t_k \to +\infty$, since $p(t_k) \ge K(t_k) > 0$, what is a contradiction. Then, the total population equalizing the susceptible one is uniformly bounded for all time. Now, assume that $N(t) \to 0$ as $t \to +\infty$. Then, it exists a strictly increasing real sequence $\{t_{1k} \in R_{0+}\}_{k \in Z_{0+}}$ such that $N(t_{1k}) \to 0$ as $R_{0+} \ni t_{1k} \to +\infty$ (so as $Z_{0+} \ni k \to +\infty$). Thus, $R_{\delta t_{1k} -} = \left\{ R_{0+} \ni \tau(\le t_k) : \left(\frac{\delta_1(\tau)N(\tau)}{p(\tau)}\right) \ge 1 \land r(t) \in R_+ \right\}$ fulfills $\lim_{k \to +\infty} \mu(R_{\delta t_{1k}-}) = +\infty$. But then $\frac{\delta_{1M} N(t_{1k})}{K(t_{1k})} \ge \frac{\delta_{1M} N(t_{1k})}{p(t_{1k})} (\to 0) \ge 1$ as $t_{1k} \to +\infty$ what is a contradiction. Then, the population cannot asymptotically extinguish.

Furthermore the above results imply that the total population is neither asymptotically strictly increasing or decreasing. Also, since $\delta_1(t)$ is not identically equal to $p(t)/N(t)$ except perhaps within some set of finite measure the total population is not everywhere constant in $R_{0+}$. Therefore, for any given $t \in R_{0+}$, it exists $t' = t'(t) > t$ such that $N(t') \ne N(t)$; i.e. either $N(t') > N(t)$ or $N(t') < N(t)$ so that the total population, equating the susceptible one for all time, exhibits a bounded non-asymptotically vanishing oscillatory behavior for any non-trivial solution of (1) in the absence of infection. □

It is interesting to consider particular models of (1) where p(t) is postulated to be a periodic function which equalizes both the total and susceptible populations when the infected one is identically zero. This reflects the particular periodic nature of the susceptible population of certain diseases versus time in the absence of infection as discussed in [26]. This is corroborated by Theorem 4.2 which proves that the total population exhibits an oscillating (rather than specifically periodic) nature in the presence of both susceptible and infected populations. Such a periodic $p: \mathbf{R}_{0+} \to \mathbf{R}_+$ is then used for the general model SIS model (1) which reflects the presence of a potential infected population by the disease.

**Assertion 4.3**. The infection-free solution of (1) has an everywhere time-differentiable periodic solution $p: \mathbf{R}_{0+} \to \mathbf{R}_{0+}$ for some period $T_p \in \mathbf{R}_+$ of the form

$$p(t+\ell T_p) = p(t) = p_0(t) + K(t) \geq K_0 + \varepsilon_0 + p_0(t) \geq \varepsilon_0 > 0$$

for some $p_0: \mathbf{R}_{0+} \to \mathbf{R}$, $K: \mathbf{R}_{0+} \to \mathbf{R}_+$, $K_0, \varepsilon_0 \in \mathbf{R}_+$; $\forall t \in \mathbf{R}_{0+}$, $\forall \ell \in \mathbf{Z}_{0+}$ if and only if $\int_t^{t+T_p} r(\tau)(1-\delta_1(\tau))d\tau = 0$; $\forall t \in \mathbf{R}_{0+}$. Such a solution might be additively decomposed in such a way that $p_0: \mathbf{R}_{0+} \to \mathbf{R}$ and $K: \mathbf{R}_{0+} \to \mathbf{R}_+$ are also everywhere time-differentiable with the same period $T_p \in \mathbf{R}_+$ as that of $p: \mathbf{R}_{0+} \to \mathbf{R}_+$.

**Proof**: Since $p_0(t) \leq K_0$ and $K(t) \geq K_0 + \varepsilon_0$ for some $K_0, \varepsilon_0 \in \mathbf{R}_+$ then $p(t) = p_0(t) + K(t) \geq K_0 + \varepsilon_0 + p_0(t) \geq \varepsilon_0 > 0$. Now fix $I \equiv 0$ and $p(t) = N(t)$ to yield from (1.a):

$$\dot{p}(t) = r(t)(1-\delta_1(t))p(t); \quad p(0) = N(0) = S(0) = N_0 = S_0 \geq \varepsilon_0 > 0$$

so that the disease-free solution $p(t) = p(0) e^{\int_0^t r(\tau)(1-\delta_1(\tau))d\tau}$ is periodic in $\mathbf{R}_{0+}$ if and only if $\exists T_p \in \mathbf{R}_+$ such that for any non-trivial solution:

$$p(t+\ell T_p) = p(t); \quad \forall t \in \mathbf{R}_{0+}, \forall \ell \in \mathbf{Z}_{0+}$$

and, equivalently, if and only if $\int_t^{t+T_p} r(\tau)(1-\delta_1(\tau))d\tau = 0$; $\forall t \in \mathbf{R}_{0+}$. The minimum $T_p \in \mathbf{R}_+$ such that $p(t+T_p) = p(t)$, $\forall t \in \mathbf{R}_{0+}$ is the period of $p: \mathbf{R}_{0+} \to \mathbf{R}_+$. Note that the periodicity of p(t) leads to the constraints:

$$K(t+\ell T_p) = K(t) + p_0(t) - p_0(t+\ell T_p) \geq \varepsilon_0 + p_0(t) - p_0(t+\ell T_p) \geq \varepsilon_0; \quad \forall t \in \mathbf{R}_{0+}, \forall \ell \in \mathbf{Z}_{0+}$$

which is achieved in a simple way in particular if $p_0(t) = p_0(t+\ell T_p)$; $\forall t \in \mathbf{R}_{0+}$, $\forall \ell \in \mathbf{Z}_{0+}$, i.e. if $p_0: \mathbf{R}_{0+} \to \mathbf{R}$ has the same period as $p: \mathbf{R}_{0+} \to \mathbf{R}_+$. □

From (2.a), (8) and $p(t) = p(0) e^{\int_0^t r(\tau)(1-\delta_1(\tau))d\tau}$, the subsequent relations between the solutions of (1) and its periodic solution $N(t) = S(t) = p(t)$ follow for initial conditions $N(0) = S(0) + I(0) = p(0) + I(0) \neq 0$ with $p(0) = S(0) = N(0)|_{I(0)=0} \neq 0$ in the absence of infected population after using the relation:

$$\int_0^t r(\tau)\left(1-\frac{\delta_1(\tau)S(\tau)}{p(\tau)}\right)d\tau = \int_0^t r(\tau)\left(1-\frac{\delta_1(\tau)}{p(\tau)}\right)d\tau + \int_0^t r(\tau)\delta_1(\tau)\left(1-\frac{S(\tau)}{p(\tau)}\right)d\tau$$

$$= p(t) + \int_0^t r(\tau)\delta_1(\tau)\left(1-\frac{S(\tau)}{p(\tau)}\right)d\tau \qquad (10)$$

so that

$$S(t) = N(t) - I(t) \;;\; I(t) = I(0)e^{\int_0^t (\beta(\tau)(N(\tau)-I(\tau))-d(\tau)-\gamma(\tau))d\tau}$$

$$N(t) = \frac{p(t)}{p(0)}(p(0)+I(0))e^{\int_0^t r(\tau)\delta_1(\tau)\left(1+\frac{I(\tau)-N(\tau)}{p(\tau)}\right)d\tau}$$

$$+\frac{I(0)}{p(0)}\int_0^t p(\tau)e^{\int_0^\tau r(\tau')\delta_1(\tau')\left(1+\frac{I(\tau')-N(\tau')}{p(\tau')}\right)d\tau'}$$

$$\left(r(t-\tau)\left[\frac{(\delta_1(t-\tau)+\delta_2(t-\tau))(N(t-\tau)-I(t-\tau))}{p(t-\tau)}-1\right]-d(t-\tau)\right)$$

$$\times e^{\int_0^{t-\tau}(\beta(\tau')(N(\tau')-I(\tau'))-d(\tau')-\gamma(\tau'))d\tau'} d\tau \;;\qquad \forall t \in \mathbf{R}_{0+}$$

The infection-free solution $p: \mathbf{R}_{0+} \to \mathbf{R}_{0+}$ can be unbounded, bounded or converge asymptotically to zero (while being or not oscillatory) instead of being periodic, under the conditions of the next result, whose proof is direct from simple inspection of the solution and then omitted:

**Assertions 4.4**. The following properties hold:
**(i)** The infection-free solution is globally asymptotically stable if and only if the subsequent integral exists:

$$\int_0^\infty r(\tau)(1-\delta_1(\tau))d\tau = -\infty$$

A sufficient condition is the existence of a finite $t_a \in \mathbf{R}_+$ such that $r(\tau)(1-\delta_1(\tau))<0$; $\forall t \in [t_a, +\infty)\setminus I_a$ where $I_a$ is a (in general, non-connected) subset of $[t_a, +\infty) \subset \mathbf{R}_{0+}$ of finite measure.

**(ii)** The infection-free solution is globally stable if

$$\int_0^t r(\tau)(1-\delta_1(\tau))d\tau < +\infty \;;\; \forall t \in \mathbf{R}_{0+}$$

A sufficient condition is the existence of a finite $t_b \in \mathbf{R}_+$ such that $r(\tau)(1-\delta_1(\tau))\leq 0$; $\forall t \in [t_b, +\infty)\setminus I_b$ where $I_b$ is a (in general, non-connected) subset of $[t_b, +\infty) \subset \mathbf{R}_{0+}$ of finite measure.

**(iii)** The infection-free solution is unstable if

$$\liminf_{t \to +\infty} \int_0^t r(\tau)(1-\delta_1(\tau))d\tau = +\infty$$

A sufficient condition is the existence of a finite $t_c \in R_+$ such that $r(\tau)(1-\delta_1(\tau))>0$; $\forall t \in [t_c, +\infty) \setminus I_c$ where $I_c$ is a ( in general, non-connected) subset of $[t_c, +\infty) \subset R_{0+}$ of finite measure. □

The positive invariance of (1) is characterized in the next result.

**Assertions 4.5**. The following properties hold:

**(i)** The domain $\Omega := \left\{ (S,I): S \geq 0, I \geq 0, N \leq \min_{t \in R_{0+}} p(t)/\delta_m \right\}$ is a positively invariant set of (1).

**(ii)** For any bounded $N(0)$, the domain

$\Omega_e(N(0)) := \{(S,I): S \geq 0, I \geq 0, N \leq N(0)\}$ is a positively invariant set of (1) provided that

$0 \leq d(t) \leq 1$, $0 \leq \delta_1(t) \leq p(t)/N(t)$, $0 \leq \delta_2(t) \leq \dfrac{p(t)(1-d(t))}{r(t)N(t)}$; $\forall t \in R_{0+}$.

**Proof**: Property (i) follows directly from $N(t) \leq p(t)/\delta_m$

$$\Rightarrow \dot{N}(t) = r(t)\left(1 - \frac{G(t)}{p(t)}\right)S(t) - d(t)I(t) \leq r(t)\left(1 - \frac{\delta_m N(t)}{p(t)}\right)S(t) - d(t)I(t)$$

$$\leq r(t)\left(1 - \frac{\delta_m N(t)}{p(t)}\right)S(t) \leq r(t)\left(1 - \frac{\delta_m N(t)}{p(t)}\right)N(t) \geq 0$$

Property (ii) follows from

$$\dot{N}(t) = r(t)\left(1 - \frac{G(t)}{p(t)}\right)S(t) - d(t)I(t) + N(t) - N(t) \leq 0$$

if and only if

$$N(t) \leq \left(1 + \left(\frac{\delta_1(t)S(t)}{p(t)} - 1\right)r(t)\right)S(t) + \left(\frac{r(t)\delta_2(t)S(t)}{p(t)} + d(t)\right)I(t) \qquad □$$

Note that if p(t) is periodic then the invariant set $\Omega$ of (1) may be defined with its defining inequality being tested only on one period. Note also that $\Omega$ is not necessarily a maximal invariant set of (1).

## 5. The impulsive model

In a animal population, three potential methods commonly used for controlling the disease, [26], are: (a) the mass pulse vaccination of the susceptible population, (b) the pulse culling or isolation or both, and (c) both (pulse culling and isolation) of the infective population, (d) a combination of vaccination and culling. Define the sequence on impulsive time instants as follows:

$$Im\, p := \{t_k \in R_{0+}: (t_{k+1} - t_k) \geq T > 0, N(t_k^+) < N(t_k), k \in SI \subset Z_+\} \qquad (11)$$

so that the epidemic model (1) is valid for $t \in R_{0+} \setminus Im\, p$ and

$$S(t_k^+) = (1-p_k)S(t_k) \;;\; I(t_k^+) = (1-q_k)I(t_k); \; \forall t_k \in Im\, p \qquad (12)$$

where $\{p_k\}_{k \in SI}$ and $\{q_k\}_{k \in SI}$ are real sequences of elements in $[0,1]$ fulfilling the constraints $0 \leq p_k q_k \leq 1$; $\forall k \in SI$. It is understood under the above notation that $t_{\bar{k}} = t_k$ what simplifies the presentation of the corresponding equations. This implies, in particular, that $p_k q_k = 0$ for some $k \in SI$ then $t_k \in Im\, p$ is not an impulsive time instant for either the susceptible or the infected. If $1 > p_k q_k > 0$ for some $k \in SI$ then $t_k \in Im\, p$ is an impulsive time instant for at least one of the susceptible and the infected populations. If $1 \geq min(p_k, q_k) > 0$ then $t_k \in Im\, p$ is an impulsive time instant for both the susceptible and the infected populations. Finally, if $p_k q_k = 1$ for some $k \in SI$ implies that $t_k \in Im\, p$ is a final impulsive time instant for both the susceptible and the infected populations which extinguishes the total population so that $N(t) = 0$; $\forall t \in [t_k, \infty)$. Note that if $t_k = +\infty$ for some finite $k = k_0 \in Z_+$, then the number of impulsive time instants is finite with $card(Im\, p) = k_0$. Note also that, since $T > 0$, the existence of a unique solution of (1) together with (11)-(12) is guaranteed for any given initial conditions. Such a solution is time- differentiable on $\bigcup_{t_k \in Im\, p}(t_k, t_{k+1})$ and bounded discontinuous for $t \in Im\, p$. It is also possible to formulate equivalently the impulsive epidemic model by defining binary real functions $p, q: R_{0+} \to [0,1]$ as follows:

$$p(t) = q(t) = 0;\ \forall t \notin Im\, p \text{ and } p(t_k) = p_k \in [0,1];\ q(t_k) = q_k \in [0,1];\ \forall t_k \in Im\, p \tag{13}$$

leading to the equivalent epidemic SIS model (1) together with (11), (13) and:

$$S(t^+) = (1 - p(t))S(t)\ ;\ I(t^+) = (1 - q(t))I(t);\ \forall t \in R_{0+} \tag{14}$$

with at least one of the susceptible or infected trajectory solutions being discontinuous at time t if and only if $t \in Im\, p$. This implies that $p(t) = 0$ if either $t \notin Im\, p$ or if $t \in Im\, p$ with $p_k = 0$; i.e. $t = t_k$ is an impulsive time instant only for the infected. In the same way, $q(t) = 0$ if either $t \notin Im\, p$ or if $t \in Im\, p$ with $q_k = 0$; i.e. $t = t_k$ is an impulsive time instant only for the susceptible. The following positivity result holds from Theorem 2.2, (11) and (13) since $0 \leq p_k q_k = p(t_k)q(t_k) \leq 1; \forall k \in SI$, $t_k \in Im\, p$ and $p(t) = 0;\ \forall t \notin Im\, p$:

**Theorem 5.1** *(Positivity of the solutions of the impulsive epidemic model)*. The solutions of the impulsive SIS epidemic model (1), (11), (13) and (14) are nonnegative for all time positive if Theorem 2.2 holds. □
The following result is direct from Theorem 2.2 and Theorem 3.1 since the susceptible, infected and total populations are non larger than the impulse- free counterparts since their numbers never increase at the impulsive time instants.

**Theorem 5.2** *(Boundedness of the solutions of the impulsive epidemic model and stability)*. If Theorem 3.1 holds then the susceptible, infected and total populations of the impulsive epidemic SIS model (1), (11) and (13)-(14) remain bounded for all time.                                                                                               □

However, it turns out that the boundedness of the total and partial populations might also be guaranteed for the impulsive model under weaker conditions than the impulsive- free counterpart by relaxing some of the positivity conditions for the parameterizing functions to be compensated by sufficient culling actions at impulsive time instants as it is now addressed in the following. Decompose the set of impulsive time instants as the (in general, non-disjoint union of the sets of impulsive time instants of the susceptible and the infected solutions of (1), i.e. $Im\, p = Im\, p(S) \cup Im\, p(I)$, where:

$$Im\, p(S) := \{t_k \in \mathbf{R}_{0+} : (t_{k+1} - t_k) \geq T > 0, S(t_k^+) < S(t_k), k \in SI \subset \mathbf{Z}_+\} \qquad (15)$$

$$Im\, p(I) := \{t_k \in \mathbf{R}_{0+} : (t_{k+1} - t_k) \geq T > 0, I(t_k^+) < I(t_k), k \in SI \subset \mathbf{Z}_+\} \qquad (16)$$

In the same way, we can define for each $t \in \mathbf{R}_{0+}$, the sets of impulsive time instants in the real intervals $[0, t)$ and $[0, t]$ as

$$Im\, p(t) := \{t_k (<t) \in \mathbf{R}_{0+} : (t_{k+1} - t_k) \geq T > 0, N(t_k^+) < N(t_k), k \in SI \subset \mathbf{Z}_+\} = Im\, p(S(t)) \cup Im\, p(I(t)) \subset Im\, p(t^+)$$

$$Im\, p(t^+) := \{t_k (\leq t) \in \mathbf{R}_{0+} : (t_{k+1} - t_k) \geq T > 0, N(t_k^+) < N(t_k), k \in SI \subset \mathbf{Z}_+\} = Im\, p(S(t^+)) \cup Im\, p(I(t^+))$$

$$Im\, p(S(t)) := \{t_k (<t) \in \mathbf{R}_{0+} : (t_{k+1} - t_k) \geq T > 0, S(t_k^+) < S(t_k), k \in SI \subset \mathbf{Z}_+\} \subset Im\, p(S(t^+))$$

$$Im\, p(S(t^+)) := \{t_k (\leq t) \in \mathbf{R}_{0+} : (t_{k+1} - t_k) \geq T > 0, S(t_k^+) < S(t_k), k \in SI \subset \mathbf{Z}_+\}$$

$$Im\, p(I(t)) := \{t_k (<t) \in \mathbf{R}_{0+} : (t_{k+1} - t_k) \geq T > 0, I(t_k^+) < I(t_k), k \in SI \subset \mathbf{Z}_+\} \subset Im\, p(I(t^+))$$

$$Im\, p(I(t^+)) := \{t_k (\leq t) \in \mathbf{R}_{0+} : (t_{k+1} - t_k) \geq T > 0, I(t_k^+) < I(t_k), k \in SI \subset \mathbf{Z}_+\} \qquad (17)$$

Note that it is not necessary to specify if the pulse culling at a time instant $t$ is for both populations or just for one of them since p(t) =0 means no culling for the susceptible and q(t)=0 means no culling for the infected and p(t)=q(t)=0 means $t \notin Im\, p$. The solution of the impulsive system is now obtained :

**Theorem 5.3** *(Weaker boundedness stability conditions of the impulsive epidemic model).* The following properties hold:

(i) All positive solutions of the impulsive SIS model under bounded initial conditions are bounded for all time for any $\delta_m > 0$.

(ii) All the positive solutions of the impulsive SIS model under bounded initial conditions are bounded for all time for any $\delta_m \geq 0$ if

$$\lim_{k \to +\infty} \sup \left( \int_0^{t_k + t} r(\tau) d\tau - \sum_{t_i \in Im\, p(t_k^+)} |\ln[1 - w(t_i)]| \right) \leq 0 \qquad \square$$